\documentclass[11pt]{amsart}

\usepackage{amsmath,amsfonts}
\usepackage{amsthm}
\usepackage{enumerate}

\usepackage[dvips]{color}

\definecolor{refkey}{gray}{0.5}
\definecolor{labelkey}{gray}{0.5}

\usepackage{graphicx}
\usepackage{caption}
\usepackage{psfrag}

\newtheorem{theorem}{Theorem}[section]
\newtheorem{proposition}[theorem]{Proposition}
\newtheorem{lemma}[theorem]{Lemma}
\newtheorem{corollary}[theorem]{Corollary}

\newtheorem*{conjecture}{Proposition}

\theoremstyle{definition}
\newtheorem{definition}[theorem]{Definition}
\newtheorem{example}[theorem]{Example}

\theoremstyle{remark}
\newtheorem{remark}[theorem]{Remark}

\numberwithin{equation}{section}

\newcommand{\nC}{\mathbb C}
\newcommand{\nD}{\mathbb D}
\newcommand{\nM}{\mathbb M}
\newcommand{\cL}{{\mathcal L}}
\newcommand{\cP}{{\mathcal P}}

\newcommand{\ds}{d_F}

 \DeclareMathOperator{\re}{Re}
\DeclareMathOperator{\aaa}{sep} 
\DeclareMathOperator{\sep1}{sep}

\begingroup
\def\2#1{\ifnum#1<10 0\fi\the#1}
\xdef\isodayandtime%
{\the\year-\2\month-\2\day\space\2{\count0}:\2{\count2}}
\xdef\now{\space\2{\count0}:\2{\count2}}
\endgroup

\begin{document}

\title[Stability of roots of polynomials]%
{Stability of roots of polynomials \\
under linear combinations of derivatives}

\author{Branko \'{C}urgus}
\address{Department of Mathematics,
Western Washington University, \newline \hspace*{3mm} Bellingham WA
98225, USA} \email{curgus@wwu.com}

\author{Vania Mascioni}
\address{Department of Mathematical Sciences, Ball State
University, \newline \hspace*{3mm} Muncie, IN 47306-0490, USA}
\email{vmascioni@bsu.edu}



\date{\today \ \ \now}

\begin{abstract}
Let $T=\alpha_0 I + \alpha_1 D + \cdots+\alpha_n D^n$, where $D$ is
the differentiation operator and $\alpha_0\not= 0$, and let $f$ be a
square-free polynomial with large minimum root separation. We prove
that the roots of $Tf$ are close to the roots of $f$ translated by
$-\alpha_1/\alpha_0$.
\end{abstract}

\maketitle

\section{Introduction}

Let $n$ be a positive integer. Denote by $\cP_n$ the
$(n+1)$-dimensional complex vector space of all polynomials of degree
at most $n$. By $\cP_0$ denote the set of all constant polynomials.
For a non-zero $\,f \in \cP_n$ let
\[
Z(f) = \bigl\{w\in\nC:f(w)=0\bigr\}
\]
be the set of the roots of $f$, where $\nC$ is the field of complex
numbers. It is useful to extend the field operations from $\nC$ onto
the nonempty subsets of $\nC$ in the following standard way:
for $A, B \subset \nC$ we set
\[
A+B = \bigl\{u+v:u\in A,v\in B\bigr\} \ \ \ \text{and} \ \ \ A\, B =
\bigl\{uv:u\in A, v\in B\bigr\}.
\]
For $r > 0$ we set $\nD(r)=\bigl\{z \in \nC:|z| \leq r\bigr\}$. For
example, the set $Z(f)+\nD(r)$ is the union of the closed disks of
radius $r$ centered at the roots of $f$.

Let $\cL(\cP_n)$ be the set of all linear operators from $\cP_n$ to
$\cP_n$.  How does an operator $T \in \cL(\cP_n)$ perturb the roots
of polynomials? To illustrate what we mean by this question consider
two simple linear operators on $\cP_n$.  Let $\alpha, t \in \nC, \,
t\neq 0$. For $f\in\cP_n$ we set
\begin{align*}
\bigl(S(\alpha) f\bigr)(z) &:= f(\alpha + z), \\
\bigl(H(t) f\bigr)(z) &:= f(z/t).
\end{align*}
Then for all non-constant $f\in\cP_n$ we clearly have
\begin{align*}
Z\bigl(S(\alpha) f\bigr) &= \{-\alpha\} + Z(f), \\
Z\bigl(H(t) f\bigr) &= \{t\} \, Z(f).
\end{align*}
Hence, for these two classes of operators there is a simple
relationship between the roots of the original polynomial and the
roots of its image. In contrast, for the differentiation operator $D:
\cP_n \to \cP_n$ and $f \in \cP_n$ there is no simple relation
between $Z(Df)$ and $Z(f)$. The classical Gauss-Lucas theorem and its
many improvements address this question; see for example the
excellent monograph on the subject \cite{RS} and \cite{MS} for a
recent development.

In this article we explore the relative position of $Z(Tf)$ in
relation to $Z(f)$ for those invertible $T \in \cL(\cP_n)$ that are
linear combinations of $I, D, \ldots, D^n$. To that end we first
define a simple measure of the perturbation of the roots under $T\in
\cL(\cP_n)$. For a non-constant polynomial $f\in\cP_n$ set
\[
R_T(f) := \min\bigl\{r > 0 : Z(Tf) \subset Z(f)+\nD(r)\bigr\}.
\]
Clearly, for the monomial $\phi_n(z) := z^n$ we have
\begin{equation*}
R_T(\phi_n) = \max\bigl\{|v| : v \in Z(T \phi_n)\bigr\}.
\end{equation*}

In \cite{CM4} we proved that the following two statements are
equivalent (see Theorem~\ref{tCA} below): \vspace*{6pt}
\begin{enumerate}[(i)]
\itemsep 6pt
\item
For all non-constant polynomials $f \in \cP_n$ we have
\[
R_T(f) \leq R_T(\phi_n).
\]
\item
There exist $\alpha_{0}, \alpha_{1}, \ldots ,\alpha_{n} \in \nC$ such
that
\begin{equation} \label{eqT}
\alpha_{0} \neq 0  \quad \text{and} \quad T = \alpha_{0} I +
\alpha_{1} D + \cdots + \alpha_{n} D^{n}.
\end{equation}
\end{enumerate}
Hence, an operator $T$ given by \eqref{eqT} perturbs the roots of
$\phi_n$ the most, as measured by $R_T$. Since $\phi_n$ only has one
root of multiplicity $n$ it is plausible to surmise that $n$ distinct
roots that are far apart from each other will be perturbed
considerably less by $T$. But is this correct?

The following few lines of {\em Mathematica} code will help explore
this question. In the code below, {\tt T} stands for a list
\[
\bigl\{ \alpha_0,\alpha_1,\ldots,\alpha_n\bigr\}, \quad \alpha_0 \neq
0,
\]
of the $n+1$ coefficients of the operator $T$ in \eqref{eqT} and {\tt
W} stands for the list of the $n$ roots of an $f \in \cP_n$. First we
calculate the roots of $Tf$ and name them {\tt TW}:

\begin{verbatim}
TW := z /. NSolve[
     T.(D[Times@@(z - W),{z, #}]&/@Range[0,Length[W]]) == 0,z
                 ]
\end{verbatim}
Then we plot the roots {\tt W} as gray points and the perturbed roots
{\tt TW} as black points by the following command:

\begin{verbatim}
 Show[Graphics[{
   {AbsolutePointSize[5], GrayLevel[0.7],
                           Point[{Re[#], Im[#]}]&/@ W},
   {AbsolutePointSize[5],  Point[{Re[#], Im[#]}]&/@ TW}
                }]]
\end{verbatim}

Applying these two commands to sets {\tt W} consisting of distinct
points that are far apart from each other the first author observed
the following surprising fact: the points of {\tt TW} were very close
to the points $\{-\alpha_1/\alpha_0\}+${\tt W}. That is, the points
of {\tt W} had been essentially translated by $-\alpha_1/\alpha_0$.
This kind of numerical experiment leads us then to believe that if
the roots of $f \in \cP_n$ are simple and far apart from each other,
then the roots of $Tf$ are close to the roots of
$S\bigl(\alpha_1/\alpha_0\bigr)f$.

Next we explore this conjecture with three simple examples. Let $a >
(n!)^{1/n}$, $\psi_{a,n}(z) := z^n - a^n$, and $T = I + D^n$. Here
$\alpha_1 = 0$, so the roots of $T\psi_{a,n}$ and $\psi_{a,n}$ should
be close for large $a$. Since $(T\psi_{a,n})(z) = z^n - (a^n-n!)$, if
we pair the roots of $\psi_{a,n}$ and $T\psi_{a,n}$ with the same
argument, then their moduli differ by $a - \bigl(a^n-n!\bigr)^{1/n}$.
This quantity indeed tends to $0$ as $a\to +\infty$. For an even $n >
2$ one can also consider $T = I + D^{n/2}$ and arrive at the same
conclusion. For $T = I + D$ we have $(T\psi_{a,n})(z) = z^n + n
z^{n-1} - a^n$ and $\alpha_1/\alpha_0 = 1$. We proceed with $n=2$,
since the roots of $T\psi_{a,2}$ are easily calculable only for this
case. Then
\[
Z\bigl(T\psi_{a,2}\bigr) = \Bigl\{-1-\sqrt{1+a^2}, -1+\sqrt{1+a^2}
\Bigr\}.
\]
To test our conjecture in this case, we consider the quantities
\begin{equation}\label{eqcl}
-a-\bigl( -1-\sqrt{1+a^2} \bigr) \quad \text{and} \quad a- \bigl(
-1+\sqrt{1+a^2}\bigr).
\end{equation}
Since both quantities converge to $1$ as $a\to +\infty$ our
conjecture is confirmed in this case as well. We offer Figures~1
and~2 as evidence that the same is true for $n=5$. In the figures the
gray points mark the roots of $\psi_{45,5}(z) = z^{5}-45^5$, the
crosses mark the points in $\{-1\}+Z(\psi_{45,5})$, while the black
dots mark the roots of $(T\psi_{45,5})(z) = z^{5} + 5z^4 -45^5$.

\noindent
\begin{figure}
\begin{minipage}[t]{.47\linewidth}
      \psfrag{p1}[][]{\begin{picture}(0,0)
            \put(-18, -10){\makebox(0,0)[l]{$z_1$}}
            \put(218,5){\vector(-1,0){190}}
                      \end{picture}}
 \psfrag{p2}[][]{\begin{picture}(0,0)
            \put(-18, 10){\makebox(0,0)[l]{$z_2$}}
                      \end{picture}}
 \psfrag{p3}[][]{\begin{picture}(0,0)
            \put(-27, -3){\makebox(0,0)[l]{$z_3$}}
            \put(110,-3){\vector(-1,0){100}}
                      \end{picture}}
 \psfrag{p4}[][]{\begin{picture}(0,0)
            \put(-27, 3){\makebox(0,0)[l]{$z_4$}}
            \put(110,2){\vector(-1,0){110}}
                      \end{picture}}
 \psfrag{p5}[][]{\begin{picture}(0,0)
            \put(-27,-3){\makebox(0,0)[l]{$z_5$}}
            \put(45,18){\vector(-3,-1){44}}
                      \end{picture}}

  \resizebox{\linewidth}{!}{\includegraphics{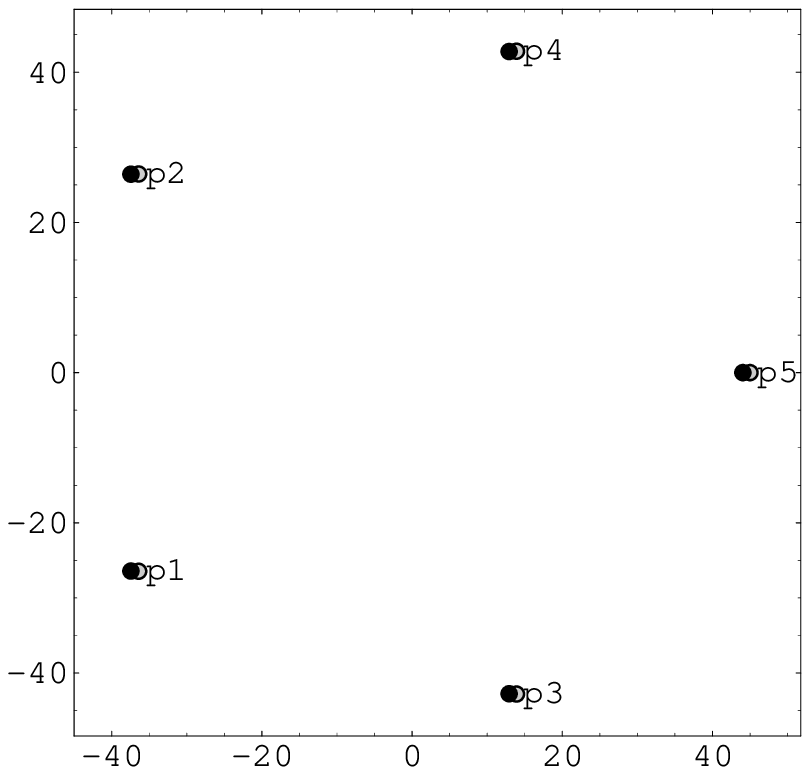}}
  \caption{$n~=~5$~and~$a~=~45$}
 \end{minipage}
    \hspace*{-5pt}  
\begin{minipage}{.46\linewidth}
\vspace*{-2.05in}
 \resizebox{\linewidth}{!}{\includegraphics{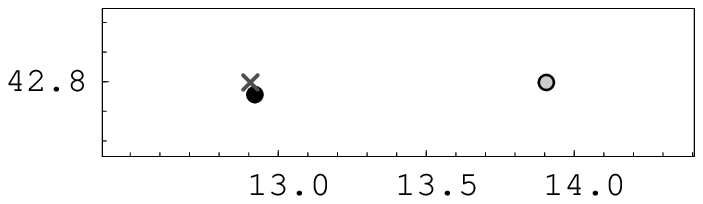}}
 \resizebox{\linewidth}{!}{\includegraphics{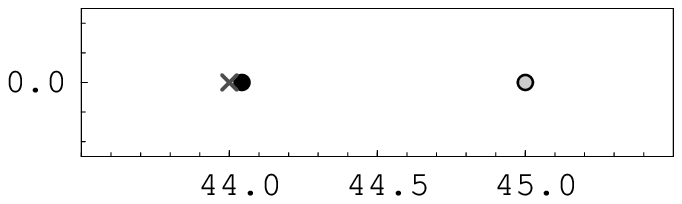}}
     \resizebox{\linewidth}{!}{\includegraphics{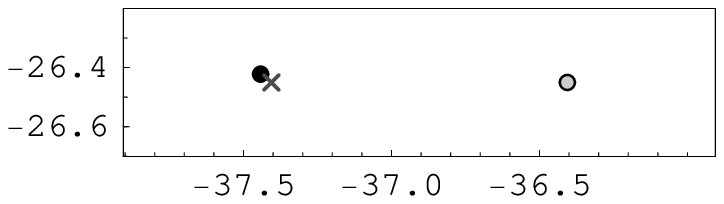}}
       \resizebox{\linewidth}{!}{\includegraphics{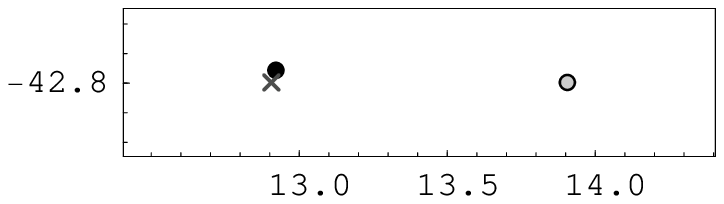}}
  \caption{A zoom-in on roots}
\end{minipage}
\end{figure}

To formulate our conjecture as a formal statement we need to define a
quantity that will play the role of $a$ in the above simple examples
and a quantity that will measure how close the roots are. To avoid
the discontinuities caused by collapsed multiple roots, we will only
consider polynomials $f$ in $\cP_n$ with simple roots. Not
surprisingly,  it turns out that the appropriate generalization for
$a$ is the minimal distance between points of $Z(f)$ and points of
$Z(f')$. We denote this quantity by $\tau(f)$, see
Definition~\ref{def1} below. As a measure how close the roots of two
polynomials with the same degree are we will use the Fr\'{e}chet
distance $d_F$: this distance is obtained by pairing the roots of two
polynomials to get the smallest possible maximal distance between the
paired roots, see Definition~\ref{def2}. Then our conjecture leads us
to the following proposition which is one of the main results of this
article:
\begin{conjecture}
Let $T$ be given by \eqref{eqT}. Then for each $\epsilon > 0$ there
exists $C_{T}(\epsilon) > 0$ such that for all $f \in \cP_n$ with
simple roots we have
\[
\tau(f) > C_{T}(\epsilon) \quad \Rightarrow \quad d_{F}\!\bigl(
Z\bigl(S\bigl(\alpha_1/\alpha_0\bigr)f\bigr), Z(T f) \bigr) <
\epsilon.
\]
\end{conjecture}

This proposition is proved as Corollary~\ref{dsTasy} in
Section~\ref{sprls}. In Section~\ref{sspr} we collect the necessary
definitions and background. Section~\ref{sRr} deals with operators $T
\in \cP_n$ given by \eqref{eqT} with $\alpha_1 = 0$. Here,
Theorem~\ref{lmt} provides the key step towards the proof of the
Proposition stated above. In Section~\ref{sprls} we prove more
results involving Fr\'{e}chet distance. We conclude with a few examples
in Section~\ref{see}.

A natural application of our results is towards a better estimate of
the regions where the roots of polynomials are located. In
Corollary~\ref{cRS} we illustrate the possible range of applications
by considering the perturbations induced by operators of the form
$T=I+\alpha D$, which have been traditionally of interest in the
study of polynomial roots, see the books \cite{Ma} and \cite{RS} for
details. This result is a variation on the often quoted classical
result by Takagi, see \cite[part (VI)]{Ta} or
\cite[Corollary~5.4.1~(iii)]{RS}. Our result, in some cases, gives
more precise information than the classical result. This is
illustrated in Example~\ref{ex3} which also gives a precise
explanation of the behavior illustrated in Figures~1 and~2 above.

\section{Definitions and preliminaries} \label{sspr}

\begin{definition} \label{def1}
Let $p$ be a polynomial of degree $n$, $n\geq 2$, and assume that $f$
has at least two distinct roots. Define
\begin{equation*} 
\sep1(f) := \min\bigl\{|w - v| \ : \ w, v \in Z(f), w \neq v \bigr\}
\end{equation*}
and
\begin{equation*} 
\tau(f) := \min\bigl\{|w - v| \ : \ w \in Z(f), v \in
Z(f')\!\setminus\!\{w\} \bigr\} \,.
\end{equation*}
\end{definition}

We will need the following inequality which was established in
\cite[Theorem~4]{CM1}.

\begin{theorem} \label{omegatau}
Let $f$ be a polynomial of degree $n$, $n\geq 2$, and assume that $f$
has at least two distinct roots. Then
\begin{equation*} 
   \frac{1}{n}\,\sep1(f) \leq \tau(f) \leq
   \frac{1}{2\sin(\pi/n)}\,\sep1(f).
\end{equation*}
\end{theorem}

\begin{remark}
The quantity $\sep1(f)$ is known as the minimum root separation of a
polynomial $f$. In \cite{CM1} we denoted it by $\omega(f)$. Finding
lower estimates for $\aaa(f)$ is important since computing time
required by an algorithm to isolate the roots of $f$ depends
inversely on $\aaa(f)$, see for example \cite{CH}. In this sense
roots of polynomials with large $\aaa(f)$ are easy to find. Our
results below indicate that such  roots are also stable under the
linear transformations of the form \eqref{eqT}.
\end{remark}

\begin{remark}
Most of our results below are formulated in terms of the quantity
$\tau(f)$. Since we consider $n$ fixed throughout,
Theorem~\ref{omegatau} yields that analogous results hold when
$\tau(f)$ is replaced by $\sep1(f)$.
\end{remark}

Next we define the Fr\'{e}chet distance (see \cite[Section~3]{CM3}
for details and references) between multisets with the same
cardinality.

\begin{definition} \label{def2}
Let $m$ be a positive integer and put $\nM = \{1,\ldots,m\}$.  By
$\Pi_m$ we denote the set of all permutations of $\nM$.  For two
functions $u, v :\nM \to \nC$ we define
\begin{equation*}
d_F(u,v) := \min_{\sigma \in \Pi_m} \max_{k\in\nM}\, \bigl|u(k) -
v(\sigma(k)) \bigr|.
\end{equation*}
The adaptation of this definition to $d_F(A,B)$, where both $A$ and
$B$ are multisets of size $n$, is then straightforward.
\end{definition}

The following theorem is a combination of Theorems~9.1 and~11.1 from
\cite{CM4}. The reader should notice that $R_T(f) = d_h(Z(f),Z(Tf))$,
$d_h$ being the notation used in \cite{CM4}.

\begin{theorem} \label{tCA}
Let $T \in \cL(\cP_n)$. The following statements are equivalent:
\begin{enumerate}[{\rm (a)}]
\itemsep 6pt
\item \label{itCA1}
There exist $\alpha_{0}, \alpha_{1}, \ldots ,\alpha_{n} \in \nC$ such
that
\begin{equation} \label{eqTa}
T = \alpha_{0} I + \alpha_{1} D + \cdots + \alpha_{n} D^{n}, \ \
\alpha_{0} \neq 0.
\end{equation}
\item \label{itCA3}
For all non-constant polynomials $f \in \cP_n$ we have
\[
   R_T(f)
   \leq R_T(\phi_n).
\]
\item \label{itCA4}
For all $f \in \cP_n$ we have $\deg(Tf) = \deg(f)$ and there exists a
constant $C_T'$ such that $\ds\bigl(Z(f),Z(Tf)\bigr) \leq C_T'$ for
all non-constant $f \in \cP_n$.
\end{enumerate}
\end{theorem}

If $T\in \cL(\cP_n)$ is given by \eqref{eqTa}, then the smallest
constant $C_T'$ which satisfies (\ref{itCA4}) in Theorem~\ref{tCA} is
denoted by $K_F(T)$.

\section{Perturbations by a special class of operators} \label{sRr}

The following theorem is the key result which is used in the rest of
the article. The main point of interest is that it treats
perturbations by operators $T$ of the type \eqref{eqTa} where the $D$
term in the expansion of $T$ is missing. In the rest of the article
we assume $n \geq 2$.

\begin{theorem} \label{lmt}
Let $\alpha_{1}, \ldots ,\alpha_{n} \in \nC$ and let $T\in
\cL(\cP_n)$ be given by
\begin{equation} \label{eqTa0}
T = I + \alpha_{1} D + \cdots + \alpha_{n} D^{n}.
\end{equation}
Then $\alpha_1 = 0$ if and only if there exists a constant $\Gamma_T
> 0$ such that for all $f \in \cP_n$ with simple roots we have
\begin{equation} \label{eqRTtb}
  {\tau(f)} \, R_T(f) < \Gamma_T.
\end{equation}
\end{theorem}
\begin{proof}
Assume that $T$ is given by \eqref{eqTa0} and that $\alpha_1=0$. Let
$f \in \cP_n$ be a polynomial with simple roots. We proceed with the
construction of $\Gamma_T$ in two steps. First, if $\tau(f) \leq 2
R_T(\phi_n) + 1$, then Theorem~\ref{tCA} immediately yields
\begin{equation} \label{eql0}
\tau(f) \leq 2\, R_T(\phi_n) + 1 \quad \Rightarrow \quad
    R_T(f)\tau(f) \leq  R_T(\phi_n) \bigl( 2R_T(\phi_n) + 1 \bigr).
\end{equation}
Next, we consider the case
\begin{equation} \label{eqla}
  \tau(f) > 2 R_T(\phi_n) + 1.
\end{equation}
Our goal is to construct an upper bound for $\tau(f) R_T(f)$ which
does not depend on $f$. By the definition of $R_T(f)$ there exist
$x_1\in Z(f)$ and $v_1\in Z(Tf)$ such that $R_T(f)=|x_1-v_1|$. In
addition, we can choose $x_1$ and $v_1$ so that
\begin{equation} \label{minxv}
  R_T(f) =  |x_1-v_1| \leq |x-v_1| \quad \text{for all} \quad x \in Z(f).
\end{equation}

The polynomial $g$ and the numbers $b_0, b_1, \ldots, b_n$ defined by
\[
g(z) : = \bigl(S(v_1)f\bigr)(z)  = b_0 + b_1 z + \frac{b_2}{2!} z +
\cdots + \frac{b_n}{n!} z^n
\]
will play an important role in the proof. Clearly $Z(g) = \{-v_1\} +
Z(f)$. Hence, by \eqref{minxv}, $w_1 = x_1 - v_1$ is a root of $g$
with smallest modulus.

Next we explore the roots of $g'$. Note that $Z(g') =
\{-v_1\}+Z(f')$. Let $z\in Z(f')$ be arbitrary. Since the roots of
$f$ are simple, $z\not= x_1$. By the definition of $\tau(f)$,
\eqref{minxv}, and \eqref{eqla} we have
\begin{equation*} 
\begin{split}
|z-v_1| & \geq  |z-x_1| - |x_1-v_1| \\
        & \geq \tau(f) - R_T(f) \\
        & >  2 R_T(\phi_n) + 1 - R_T(\phi_n) \\
        & = R_T(\phi_n)+1.
\end{split}
\end{equation*}
Thus,
\begin{equation} \label{eqzTffp1}
|u| > R_T(\phi_n) + 1 > 0 \ \ \ \text{for all} \ \ \  u \in Z(g').
\end{equation}
In particular,
\begin{equation*} 
b_1 = g'(0) \neq 0.
\end{equation*}

Now recall that $v_1 \in Z(T f)$, $\alpha_1=0$, and observe that
$S(v_1)$ and $T$ commute, to deduce
\begin{equation} \label{eqTf1}
0 = (Tf)(v_1) =  \bigl(S(v_1)T f \bigr)(0) = \bigl(T g \bigr)(0) =
b_0 + \sum_{k=2}^n \, \alpha_k b_k.
\end{equation}
Since $w_1$ is a root of $g$ with smallest modulus, Vi\`ete's
formulas for $g$ imply
\begin{equation*}
   |b_1| |w_1| \leq n |b_0|.
\end{equation*}
Together with \eqref{eqTf1} this yields
\begin{equation} \label{nse}
   |w_1| \leq  \frac{n}{|b_1|} \sum_{k=2}^n \, {|\alpha_k| |b_k|} .
\end{equation}

Next we consider the polynomial
\[
h(z):= z^{n-1}g'(1/z) = b_1 z^{n-1} + b_2 z^{n-2} + \ldots +
\frac{b_{n-1}}{(n-2)!} z + \frac{b_n}{(n-1)!}.
\]
We recall that $b_1 \not= 0$ and we notice that for $k = 2, \ldots,
n$, the number $|b_k|/\bigl(|b_1|(k-1)!\bigr)$ is the modulus of the
coefficient of $z^{n-k}$ in the monic polynomial $h/b_1$. Further,
let $u_1$ be a root of $g'$ with minimal modulus. Then the largest
modulus among the roots of $h$ is $1/|u_1|$ and, again by Vi\`ete's
formulas, we have
\begin{equation} \label{eqmins}
   \left|\frac{b_k}{b_1(k-1)!}\right| \leq \binom{n-1}{k-1}
   \frac{1}{|u_1|^{k-1}},  \ \ \ k = 2, \ldots, n .
\end{equation}

Since the function
 $x \mapsto {\bigl((1+x)^{n-1} - 1\bigr)}/{x}, x > 0$,
is increasing, \eqref{eqzTffp1} implies
\begin{equation} \label{eqmi}
|u_1| \left(\!\! \left(1+\frac{1}{|u_1|}\right)^{n-1}\!\!-1\! \right)
\leq
 \bigl( R_T(\phi_n) + 1 \bigr)\!
 \left(\!\! \left(\frac{R_T(\phi_n) + 2}{R_T(\phi_n) + 1}\right)^{n-1}\!\! - 1\!\! \right).
\end{equation}

Next we use \eqref{nse}, \eqref{eqmins} and  \eqref{eqmi} to
establish an upper estimate for $ |u_1w_1| $:
\begin{align*}
|u_1 w_1| &\leq
   |u_1| \, n  \sum_{k=2}^n \,\frac{|\alpha_k| |b_k|}{|b_1|}\\
        &\leq
    |u_1| \,  n \left(\max_{2\leq k \leq n}
        {|\alpha_k|(k-1)!}\right)
        \sum_{k=2}^n \left|\frac{b_k}{b_1(k-1)!}\right| \\
        &\leq
    |u_1| \, n \left(\max_{2\leq k \leq n}
        {|\alpha_k|(k-1)!}\right) \sum_{k=1}^{n-1}
        \binom{n-1}{k}\frac{1}{|u_1|^{k}}\\
        & =
    |u_1| \, n \left(\max_{2\leq k \leq n}
        {|\alpha_k|(k-1)!}\right)
        \left(\!\!\left(1+\frac{1}{|u_1|}\right)^{n-1}-1\!\right) \\
    & \leq n \left(\max_{2\leq k \leq n}
        {|\alpha_k|(k-1)!}\right)
         \bigl( R_T(\phi_n) + 1 \bigr)\!
 \left(\!\! \left(\frac{R_T(\phi_n)+2}{R_T(\phi_n) + 1}\right)^{n-1}\!\! - 1\!\! \right).
\end{align*}
For further reference we denote the last expression by $\Gamma'_T$.
Hence we proved
\begin{equation} \label{eqwuLE}
|u_1 w_1| \leq \Gamma'_T.
\end{equation}

Recall that  $u_1\in Z(g')$ and $w_1\in Z(g)$. Since the roots of $g$
are simple, $w_1\not= u_1$ and hence
\[
|u_1-w_1| \geq  \tau(g) = \tau(f).
\]
By Theorem~\ref{tCA} and \eqref{eqla}
\[
R_T(f) \leq R_T(\phi_n) < \tau(f)/2.
\]
Now, the triangle inequality and \eqref{minxv} yield
\begin{align*}
    |u_1| & \geq |u_1-w_1| - |w_1| \\
     & \geq \tau(f) - R_T(\phi_n)\\
     & > \tau(f)/2 .
\end{align*}
Since by our choice $|w_1|=|x_1-v_1|=R_T(f)$, \eqref{eqwuLE} now
yields
\begin{equation} \label{mmm}
  \tau(f) R_T(f)  < 2\, |u_1| \, |w_1| \leq
  2 \, \Gamma'_T.
\end{equation}
Thus, we have proved the implication \eqref{eqla} $\Rightarrow$
\eqref{mmm}. To conclude the proof of the ``only if'' part of the theorem, set
\begin{equation*} 
\Gamma_T  :=
 \max \bigl\{ R_T(\phi_n) \bigl(2R_T(\phi_n)+1\bigr),  2\, \Gamma'_T \bigr\}.
\end{equation*}
With this $\Gamma_T$ the implication \eqref{eqla} $\Rightarrow$
\eqref{eqRTtb} clearly holds. Finally, recall the implication
\eqref{eql0} and the ``only if'' part is proved.

To prove the ``if'' part of the theorem, assume that $T$ is given by
\eqref{eqTa0} and that \eqref{eqRTtb} holds for all $f\in\cP_n$ with
simple roots. Let $a > 0$ be arbitrary and as before
$\psi_{a,2}(z)=z^2 - a^2$. Then $a$ and $-a$ are the roots of
$\psi_{a,2}$ and $\tau(\psi_{a,2}) = a$. Also,
\[
(T\psi_{a,2})(z) = z^2 + 2 \alpha_1 z - a^2 + 2 \alpha_2.
\]
The roots of $T\psi_{a,2}$ are
\[
z_{a,1} = -\alpha_1 - \sqrt{\alpha_1^2 + a^2 - 2 \alpha_2} \quad
\text{and} \quad z_{a,2} = -\alpha_1 + \sqrt{\alpha_1^2 + a^2 - 2
\alpha_2 }.
\]
Since clearly
\begin{align*}
\lim_{a\to+\infty}
 \Bigl|-a - \Bigl(-\alpha_1 - \sqrt{\alpha_1^2 + a^2 - 2 \alpha_2 } \Bigr) \Bigr|
  & = |\alpha_1|, \\
\lim_{a\to+\infty}
 \Bigl|a - \Bigl(-\alpha_1 + \sqrt{\alpha_1^2 + a^2 - 2 \alpha_2 } \Bigr) \Bigr|
 & = |\alpha_1|, \\
\lim_{a\to+\infty}
 \Bigl|-a - \Bigl(-\alpha_1 + \sqrt{\alpha_1^2 + a^2 - 2 \alpha_2} \Bigr) \Bigr|
  & = +\infty, \\
\lim_{a\to+\infty}
 \Bigl|a - \Bigl(-\alpha_1 - \sqrt{\alpha_1^2 +a^2 - 2 \alpha_2} \Bigr) \Bigr|
  & = +\infty,
\end{align*}
we conclude
\[
\lim_{a\to+\infty} R_T(\psi_{a,2}) = |\alpha_1|.
\]
Since by \eqref{eqRTtb} for all $a > 0$ we have $R_T(\psi_{a,2}) \leq
\Gamma_T/a$, letting $a \to +\infty$ leads to $ |\alpha_1| = 0$. This
completes the proof.
\end{proof}

\begin{corollary} \label{clmt}
Let $\alpha_2, \ldots, \alpha_n \in \nC$.  Let $T\in \cL(\cP_n)$ be
given by
 \begin{equation*}  
T = I  + \alpha_{2} D^2 + \cdots + \alpha_{n} D^{n}.
\end{equation*}
Let $f \in \cP_n$ be a polynomial with simple roots such that
$\tau(f) > 2 R_T(\phi_n) + 1$. Then
\[
Z(Tf) \subset Z(f) + \nD\bigl(\Gamma_T'/\tau(f)\bigr),
\]
with $\Gamma_T'$ as defined in the sentence preceding \eqref{eqwuLE}.
\end{corollary}
\begin{proof}
The corollary is in fact a restatement of the implication
\eqref{eqla}$\Rightarrow$\eqref{mmm} which is proved as a part of the
proof of Theorem~\ref{lmt}.
\end{proof}

The following corollary is inspired by
\cite[Corollary~5.4.1(iii)]{RS}.

\begin{corollary} \label{cRS}
Let $\alpha \in \nC$ and let $T\in \cL(\cP_n)$ be given by $T = I
+\alpha D$. Let $f \in \cP_n$, $n\geq 2$, be a polynomial with simple
roots such that $\tau(f) > 2 |\alpha|(n-1)+1$. Then
\begin{equation} \label{eqinc}
Z(Tf) \subset \{-\alpha\} + Z(f) + \nD\bigl(\gamma_\alpha /\tau(f)
\bigr),
\end{equation}
where $\gamma_{\alpha}$ is given by
\[
2n\bigl(|\alpha|(n-1)+1\bigr)\!\!
  \left( \!\!\!\biggl(\!\frac{|\alpha|(n-1)+2}{|\alpha|(n-1)+1}\!\biggr)^{\!n-1}\!\!\! -1\!\! \right) \! %
   \max\!\left\{\!|\alpha|^k  \Bigl(\! 1- \frac{1}{k}\! \Bigr) ,
   k=2,\ldots,n\!
\right\}\!.
\]
\end{corollary}
\begin{proof}
Corollary~\ref{clmt} does not apply to the operator $T\in
\cL(\cP_n)$. Therefore we consider the composition $V =
S(-\alpha)T \in \cL(\cP_n)$ where $S(-\alpha) \in \cL(\cP_n)$ is
defined by $\bigl(S(-\alpha)f\bigr)(z) = f(z-\alpha)$. The Taylor
formula at $z$ implies that $S(-\alpha) = \sum_{k=0}^n (-\alpha)^k
D^k/k!$. Hence $S(-\alpha)$ and $T$ commute. Since $D \in \cL(\cP_n)$
is the differentiation operator on $\cP_n$ we have $D^{n+1} =0$.
These observations lead to the following expression for $V$ as a
linear combination of derivatives:
\begin{align*}
 V & =     (I + \alpha D) \sum_{k=0}^n \frac{(-\alpha)^k}{k!} D^k \\
 & = \sum_{k=0}^n \frac{(-\alpha)^k}{k!}D^k - \sum_{k=0}^n \frac{(-\alpha)^{k+1}}{k!}D^{k+1} \\
  & = I + \sum_{k=1}^n \frac{(-\alpha)^k}{k!}(1-k)D^k\\
   & = I - \frac{\alpha^2}{2} D^2 +  \frac{\alpha^3}{3} D^3 - \cdots
     + \frac{(-1)^{n+1}\alpha^n}{(n-2)!\,n}D^n.
\end{align*}
Hence, Corollary~\ref{clmt} applies to $V$. To show that $\Gamma_V' =
\gamma_\alpha$ we first calculate $(T\phi_n)(z) = z^n+\alpha n
z^{n-1}$ and deduce $Z(T\phi_n) = \{0,-\alpha n\}$. Thus,
$R_T(\phi_n) = |\alpha| n$. Also,
\[
Z(V\phi_n) = Z\bigl(S(-\alpha)T \phi_n \bigr) = \{\alpha\} +
\{0,-\alpha n\} = \{\alpha, \alpha (1-n) \},
\]
and therefore, $R_V(\phi_n) = |\alpha|(n-1)$. As we calculated the
coefficients of $V$ to be
\[
\alpha_k= \frac{(-1)^{k+1}\alpha^k}{(k-2)! \, k}, \qquad k =
2,\ldots,n,
\]
we have $\gamma_\alpha = \Gamma_V'$.

Since we assume that $\tau(f) > 2 |\alpha|(n-1)+1 = 2 R_V(\phi_n)+1$,
Corollary~\ref{clmt} yields
\[
Z(Vf) \subset Z(f) + \nD\bigl( \gamma_\alpha /\tau(f) \bigr).
\]
To obtain the inclusion
in the corollary, substitute $f$ with $S(\alpha)f$ and notice that
$\tau\bigl(S(\alpha)f\bigr) = \tau(f)$, $VS(\alpha)f =Tf$ and
$Z\bigl(S(\alpha)f\bigr) = \{-\alpha\}+Z(f)$.
\end{proof}

\begin{remark}
In our notation the inclusion proved in
\cite[Corollary~5.4.1(iii)]{RS} reads
\begin{equation} \label{eqincRS}
Z(Tf) \subset \{-\alpha\,n/2\} + Z(f) + \nD\bigl( |\alpha|\,n/2
\bigr).
\end{equation}
Notice that if $n>1$ and
\[
\tau(f) > \max\left\{ 2 |\alpha|(n-1)+1,
\frac{\gamma_\alpha}{|\alpha|(n-1)} \right\},
\]
then the intersection of the right-hand sides of \eqref{eqinc} and
\eqref{eqincRS} provides improved information about the location of the
roots of $Tf$. Moreover, if
\[
\tau(f) > \max\left\{ 2 |\alpha|(n-1)+1,
\frac{\gamma_\alpha}{|\alpha|} \right\},
\]
then \eqref{eqinc} is an improvement of \eqref{eqincRS}, the
improvement being considerable for large $\tau(f)$.
\end{remark}

\begin{example}
In this example we give a hint of the problems that can arise if we consider
polynomials with multiple roots. For $a > 0$ set
\[
   g_a(z) = z^2(z-a)^2.
\]
Let
\[
   T = I + D^2.
\]
Since $g_a'(z)=2 z (z-a) (2 z-a)$, we have $\tau(g_a) = a/2$.  On the
other hand, we have
\begin{align*}
   (Tg_a)(z) & = g_a(z) + g_a''(z) \\
       & = z^4-2 a z^3+\left(a^2+12\right) z^2-12 a z+2 a^2 \\
       & = \frac{a^4}{16}-a^2+\left(12-\frac{a^2}{2}\right)
       \left(z - \frac{a}{2} \right)^2 +\left(z - \frac{a}{2} \right)^4.
\end{align*}
Thus the roots of $Tg_a$ are symmetric with respect to the real axis
and to the vertical line $\re(z) = a/2$. So, it is sufficient to
calculate one root of $Tg_a$ with a positive imaginary part and real
part less than $a/2$. For $a > 3 \sqrt{2}$ such a root is given by
\begin{equation*}
 z_{a,1} = \frac{1}{2}
  \left( a - \sqrt{a^2- 24 - 4  i \sqrt{2 a^2 - 36}}\right),
\end{equation*}
and the other three roots are
\[
z_{a,2} = \overline{z_{a,1}}, \quad z_{a,3} = z_{a,1} + a, \quad
z_{a,4} = \overline{z_{a,1}} + a.
\]
By rationalizing and then simplifying the expression for $z_{a,1}$ it
is not difficult to prove that \[
 \lim_{a \to +\infty}  z_{a,1} = \sqrt{2}\, i .
\]
Consequently,
\begin{align*}
 \lim_{a \to +\infty} \bigl( z_{a,2} -(-\sqrt{2}\, i) \bigr) & = 0, \\
 \lim_{a \to +\infty} \bigl( z_{a,3} - ( a + \sqrt{2}\, i ) \bigr)  & =
 0,\\
 \lim_{a \to +\infty} \bigl( z_{a,4} - ( a - \sqrt{2}\, i ) \bigr)  & = 0.
\end{align*}
Since $Z(g_a) = \{0,a\}$, the last four equalities imply
\[
\lim_{a\to+\infty} R_T(g_a) = \sqrt{2}.
\]
As $\tau(g_a) = a/2$, the function $a \mapsto \tau(g_a) R_T(g_a), \,
a > 0$, is unbounded. Hence, the assumption of the simplicity of the
roots in Theorem~\ref{lmt} cannot be dropped.
\end{example}

\section{Results involving the Fr\'{e}chet distance} \label{sprls}

In the previous section we used $R_T(f)$ as a measure of the distance
between the roots of $f$ and the roots of $Tf$. In this section we
work with the Fr\'{e}chet distance between the roots of polynomials
$f$ and $Tf$. Here $Z(f)$ and $Z(Tf)$ are considered as multisets of roots.  The Fr\'{e}chet distance is defined only for nonconstant
polynomials with equal degrees. In each of the cases below, the fact
that the degrees of $f$ and $Tf$ are equal follows from
Theorem~\ref{tCA}. In particular, $S(\alpha)$ does not change the
degree of a polynomial, and thus all the Fr\'{e}chet distances used
below are well defined. Recall that the number $K_F(T)$ is defined
immediately after Theorem~\ref{tCA}. The following lemma provides a connection
with the results from Section~\ref{sRr}.

\begin{lemma} \label{lFd}
Let $T\in \cL(\cP_n)$ be given by
\begin{equation} \label{eqT1}
T = \alpha_0 I + \alpha_{1} D + \cdots + \alpha_{n} D^{n}, \quad
\alpha_{0},  \alpha_{1}, \ldots ,\alpha_{n} \in \nC, \quad \alpha_0
\neq 0.
\end{equation}
Let $f \in \cP_n$ be a polynomial with simple roots such that
\begin{equation} \label{eqcFd}
R_T(f) < 1   \quad  \text{and}  \quad  \tau(f) >
\frac{1+K_F(T)}{\sin(\pi/n)}.
\end{equation}
Then
\begin{equation} \label{eqFd=R}
\ds\bigl(Z(f),Z(Tf)\bigr) = R_T(f) .
\end{equation}
\end{lemma}
\begin{proof}
Let $f \in \cP_n$ be a polynomial with simple roots satisfying
\eqref{eqcFd}. Then Theorem~\ref{omegatau} yields
\begin{equation} \label{eqmt12}
\sep1(f) > 2 + 2 K_F(T) \geq 2.
\end{equation}
Further, by the definition of $K_F(T)$, we have
\[
d_F\bigl(Z(f),Z(Tf)\bigr) \leq K_F(T).
\]
Now assume that the degree of  $f$ is $m$ and let $Z(f) =
\{z_1,\ldots,z_m\}$, where $z_1,\ldots,z_m$ are distinct complex
numbers.  By the definition of the Fr\'{e}chet distance
$d_F\bigl(Z(f),Z(Tf)\bigr)$ the roots $w_1,\ldots,w_m$ of $Tf$
(recall that also $\deg Tf = m$) can be indexed in such a way that
\begin{equation} \label{eqmt13}
  |z_j - w_j| \leq K_F(T), \quad j = 1,\ldots,m.
\end{equation}
Let $j$ and $k$ be distinct numbers from $\{1,\ldots,m\}$. Then, $z_j
\neq z_k$ and the triangle inequality, \eqref{eqmt12} and
\eqref{eqmt13} yield
\begin{eqnarray*}
  |w_j-w_k| & = & | (w_j -z_j) + (z_j - z_k) + (z_k  - w_k) |  \\
           &\geq & |z_j-z_k| - |(w_j -z_j) + (z_k  - w_k)| \\
           &\geq & |z_j-z_k| - (|w_j-z_j| + |z_k-w_k|) \\
           &\geq & |z_j-z_k| - |z_j-w_j| - |z_k-w_k| \\
           &\geq &   \sep1(f) - 2K_F(T) \\
           & > &  2.
\end{eqnarray*}
Consequently,
\[
\sep1(Tf) > 2.
\]
Therefore, the disks $\{w_j\}+\nD(1)$ (for $j \in \{1,\ldots,m\}$)
are pairwise disjoint. Since $R_T(f) < 1$ and $\sep1(f) > 2$, in each
disk $\{w_j\}+\nD(1)$, $j \in \{1,\ldots,m\}$, there is exactly one
root of $f$. Renumber the roots of $f$ in such a way that
\[
z_j \in \{w_j\}+\nD(1), \quad j \in \{1,\ldots,m\}.
\]
Then
\[
\ds\bigl(Z(f),Z(Tf)\bigr) \leq \max \bigl\{|w_j-z_j| : j \in
\{1,\ldots,m\}\bigr\} = R_T(f) .
\]
Since clearly $R_T(f) \leq \ds\bigl(Z(f),Z(Tf)\bigr)$, \eqref{eqFd=R}
follows.
\end{proof}

\begin{theorem} \label{mainrevised}
Let $n \geq 2$ and let $T\in \cL(\cP_n)$ be given by
\begin{equation} \label{eqT0}
T = I + \alpha_{2} D^2 + \cdots + \alpha_{n} D^{n}, \quad \alpha_{2},
\ldots ,\alpha_{n} \in \nC.
\end{equation}
Then for every $\varepsilon > 0$ there exists $C_T(\varepsilon) > 0$
such that for all $f \in \cP_n$ with simple roots the following
implication holds:
\[
\tau(f) > C_T(\varepsilon) \quad \Rightarrow \quad
d_F\bigl(Z(f),Z(Tf)\bigr) < \varepsilon.
\]
\end{theorem}

\begin{proof}
Let $\Gamma_T' > 0$ be the constant introduced immediately above inequality \eqref{eqwuLE}.
For $\varepsilon > 0$ set
\begin{equation} \label{eqCe}
C_T(\varepsilon) = \max \left\{ R_T(\phi_n) + 1, \Gamma_T',
\frac{\Gamma_T'}{\varepsilon}, \frac{1+K_F(T)}{\sin(\pi/n)} \right\}.
\end{equation}

Let $f \in \cP_n$ be a polynomial with simple roots such that
\begin{equation*}  
\tau(f) > C_T(\varepsilon).
\end{equation*}
Since $\tau(f) > R_T(\phi_n)+1$, Corollary~\ref{clmt} implies
\begin{equation} \label{eqmt11}
   R_T(f) < \frac{\Gamma_T'}{\tau(f)}.
\end{equation}
Since $\tau(f) > \Gamma_T'$, $R_T(f) < 1$. Hence, Lemma~\ref{lFd}
yields
\[
\ds\bigl(Z(f),Z(Tf)\bigr) = R_T(f).
\]
Since $\tau(f) > \Gamma_T'/\varepsilon$, \eqref{eqmt11} implies
\[
R_T(f) < \varepsilon.
\]
The last two displayed relations prove the theorem.
\end{proof}

The following proposition is proved by combining the methods of
proofs of Theorems~\ref{lmt} and~\ref{mainrevised}. In the
proposition we use $\Gamma_T$ defined in Theorem~\ref{lmt}.

\begin{proposition} \label{pub}
If $T\in \cL(\cP_n)$ is given by \eqref{eqT0}, then for all $f \in
\cP_n$ with simple roots we have
\begin{equation*} 
\tau(f) \, d_F\bigl(Z(f),Z(Tf)\bigr) \leq \max
\left\{\Gamma_T, K_F(T) \Gamma_T ,  K_F(T)\frac{1 + K_F(T)}{\sin(\pi/n)}
\right\}.
\end{equation*}
\end {proposition}
\begin{proof}
Set
\[
\Gamma_T^{\prime\prime} := \max \left\{\Gamma_T, \frac{1 +
K_F(T)}{\sin(\pi/n)} \right\}.
\]
Let $f \in \cP_n$ be a polynomial with simple roots. As in the proof
of Theorem~\ref{lmt} we proceed in two steps. First, if $\tau(f) \leq
\Gamma_T^{\prime\prime}$, then, by the definition of $K_F(T)$,
\begin{equation*} 
\tau(f) \, d_F\bigl(Z(f),Z(Tf)\bigr) \leq  \Gamma_T^{\prime\prime} \,
K_F(T) \leq   \max \left\{\Gamma_T, \Gamma_T^{\prime\prime} K_F(T)
\right\}.
\end{equation*}
For the second step assume $\tau(f) >\Gamma_T^{\prime\prime}$. Recall
that by Theorem~\ref{lmt}
\begin{equation} \label{eqmt112}
  \tau(f) \, R_T(f) \leq \Gamma_T \leq \Gamma_T^{\prime\prime}.
\end{equation}
Consequently, $R_T(f) < 1$. Since also $\tau(f) >
\bigl(1+K_F(T)\bigr)/(\sin(\pi/n))$, Lemma~\ref{lFd} implies
\[
\ds\bigl(Z(f),Z(Tf)\bigr) = R_T(f).
\]
Substituting this identity in \eqref{eqmt112}, we conclude
\begin{equation*} 
\tau(f) \, d_F\bigl(Z(f),Z(Tf)\bigr) \leq \Gamma_T
 \leq \max \left\{\Gamma_T,
\Gamma_T^{\prime\prime} K_F(T) \right\},
\end{equation*}
and the proposition is proved.
\end{proof}

Theorem~\ref{mainrevised} is of course a consequence of
Proposition~\ref{pub}, but $C_T(\epsilon)$ given by \eqref{eqCe} is
smaller than the corresponding constant deduced from
Proposition~\ref{pub}.

Now we have all the tools to prove the proposition stated in the
Introduction.

\begin{corollary} \label{dsTasy}
Let $T\in \cL(\cP_n)$ be given by \eqref{eqT1}. Then for every
$\varepsilon > 0$ there exists $C_T(\varepsilon) > 0$ such that for
all $f \in \cP_n$ with simple roots we have
\[
\tau(f) > C_T(\varepsilon) \quad \Rightarrow \quad \ds\bigl(
Z(S(\alpha_1/\alpha_0)f),Z(Tf) \bigr)  < \varepsilon.
\]
\end{corollary}

\begin{proof}
To prove the corollary we observe the identity
\[
\ds\bigl( Z(S(\alpha_1/\alpha_0)f),Z(Tf) \bigr) = \ds\bigl(
Z(f),Z\bigl(S(-\alpha_1/\alpha_0)Tf\bigr) \bigr), \ \ f \in
\cP_n\setminus\cP_0.
\]
Next we notice that Theorem~\ref{mainrevised} applies to the operator
\[
T_1 :=\frac{1}{\alpha_0} S(-\alpha_1/\alpha_0)T.
\]
Then the corollary follows from Theorem~\ref{mainrevised}.
\end{proof}

A simple calculation shows that for $t > 0$ we have
\begin{equation*}
\tau\bigl(H(t)f\bigr) = t\,\tau(f),
\end{equation*}
and thus we note the following corollary.
\begin{corollary}
Let $T\in \cL(\cP_n)$ be given by \eqref{eqT1}. For an arbitrary $f
\in \cP_n$ with simple roots we have
\[
\lim_{t\to +\infty} d_{_F}\!\Bigl( Z\bigl(S(\alpha_1/\alpha_0) H(t)
f\bigl),Z\bigl(T H(t) f\bigr) \Bigr) = 0
\]
\end{corollary}

Now we return to a simple motivating example from the Introduction. Finding the limits of the expressions in \eqref{eqcl} is a standard calculus exercise. Since the expressions in \eqref{eqcl} represent distances between roots of two quadratic polynomials, the following corollary can be seen as a generalization of this standard calculus exercise.

\begin{corollary} \label{frstcor}
Let $a_1, \ldots, a_{n-1}$ be arbitrary, but fixed, complex numbers.
For $a \in \nC$ set
\[
g_a(z) = z^n + a^n \quad \text{and} \quad f_a(z) = z^n + a_{n-1}
z^{n-1} + \cdots + a_1 z + a^n.
\]
Then
\begin{equation*}
\lim_{|a|\to +\infty} d_{_F}\!\Bigl(\! Z\bigl(S(a_{n-1}/n)
g_a\bigl),Z\bigl(f_a\bigr)\!\Bigr) = 0.
\end{equation*}
\end{corollary}
\begin{proof}
Setting
\[
\alpha_0 = 1, \ \ \alpha_{k} = \frac{(n-k)!}{n!}\, a_{n-k}, \ \ k =
1,\ldots, n-1, \ \ \alpha_n = 0,
\]
in \eqref{eqT1} we have $Tg_a = f_a$. Since clearly $\tau(g_a) =
|a|$, the corollary follows from Corollary~\ref{dsTasy}.
\end{proof}

Corollary \ref{frstcor} indicates that for a monic polynomial $f_a$ with constant term $a^n$ of large magnitude the approximate location of its roots depends only on $a$, $n$, and the coefficient $a_{n-1}$. The following corollary formalizes this observation.

\begin{corollary} \label{scndcor}
Let $a_1, \ldots, a_{n-1}$ be arbitrary, but fixed, complex numbers.
For $a \in \nC$ set
\[
h_a(z) = z^n + a_{n-1} z^{n-1} + a^n \quad \text{and} \quad f_a(z) =
z^n + a_{n-1} z^{n-1} + \cdots + a_1 z + a^n.
\]
Then
\begin{equation*}
\lim_{|a|\to +\infty} d_{_F}\!\bigl( Z(h_a),Z(f_a) \bigr) = 0.
\end{equation*}
\end{corollary}
\begin{proof}
Let $g_a(z) = z^n - a^n$. By Corollary~\ref{frstcor}
\begin{align*}
\lim_{|a|\to +\infty} d_{_F}\!\Bigl(\! Z\bigl(S(a_{n-1}/n)
g_a\bigl),Z\bigl(f_a\bigr)\!\Bigr) & = 0 \\
 \intertext{and}
\lim_{|a|\to +\infty} d_{_F}\!\Bigl(\! Z\bigl(S(a_{n-1}/n)
g_a\bigl),Z\bigl(h_a\bigr)\!\Bigr) & = 0.
\end{align*}
Since by \cite[Proposition~3.1]{CM3}, $d_F$ is a metric on the space of root sets $Z(\cdot)$,
\begin{multline*}
d_{_F}\!\bigl( Z(h_a),Z(f_a) \bigr) \\
\leq d_{_F}\!\Bigl(\! Z\bigl(S(a_{n-1}/n)
g_a\bigl),Z\bigl(f_a\bigr)\!\Bigr) + d_{_F}\!\Bigl(\! Z\bigl(S(a_{n-1}/n)
g_a\bigl),Z\bigl(h_a\bigr)\!\Bigr).
\end{multline*}
The corollary follows from the last three displayed relations.
\end{proof}

\section{Estimates and examples} \label{see}

\begin{example} \label{ex1}
A simple way to see our results in action, in particular
Theorem~\ref{lmt}, is to look at an example for which we can exactly
calculate $\tau(f)$ and $R_T(f)$.  In the Introduction we considered
$T = I + D^n$. Here $\alpha_1 = 0$, so  Theorem~\ref{lmt} applies.
Further, $R_T(\phi_n)=(n!)^{1/n}$ and, for $n \geq
3$,
\begin{equation*} 
\Gamma_T  = 2 n! \bigl((n!)^{1/n} + 1\bigr) \left(\!\!
\left(\frac{(n!)^{1/n}+2}{(n!)^{1/n} + 1}\right)^{n-1}\!\! - 1\!\!
\right).
\end{equation*}
A rough estimate yields
\[
\Gamma_T < 2\bigl(\sqrt{2}+1\bigr)\, n! \left(\frac{3}{2}\right)^{n-1}  .
\]
Hence, by Theorem~\ref{lmt},
\begin{equation} \label{eqsc}
\tau(f) \, R_T(f) < 2\bigl(\sqrt{2}+1\bigr)\, n!\left(\frac{3}{2}\right)^{n-1}
\end{equation}
for all $f \in \cP_n$ with simple zeros. Let $a>0$ and $\psi_{a,n}(z)
= z^n - a^n$. Then $\tau(\psi_{a,n}) = a$ and $(T\psi_{a,n})(z) = z^n
- (a^n-n!)$. It is not difficult to see that
\[
R_T(\psi_{a,n}) = \bigl| a - \bigl(a^n - n!\bigr)^{1/n} \bigr|,
\]
where, for $0 < a < (n!)^{1/n}$, the formula $\bigl(a^n -
n!\bigr)^{1/n}$ denotes the root with argument $\pi/(2n)$. Further,
elementary considerations yield
\[
\tau(\psi_{a,n}) \, R_T(\psi_{a,n}) = a\, \bigl| a - \bigl(a^n -
n!\bigr)^{1/n} \bigr| \leq (n!)^{2/n}.
\]
This is certainly much better than \eqref{eqsc}. However,
\eqref{eqsc} holds for \textit{all} $f \in \cP_n$ with simple
zeros.
\end{example}

\begin{remark}
Example~\ref{ex1} indicates that the constant $\Gamma_T$ in
Theorem~\ref{lmt} might not be close to optimal, at all. We note that
Theorem~\ref{lmt} together with Theorem~\ref{omegatau} yields
information about the maximum root separation $\sep1(f)$:
\[
\sep1(f) R_T(f) \leq n \Gamma_T.
\]
In \cite{Co}, Collins provides us with an impressive amount of
numerical evidence for what should be some general ``ideal'' lower
bound for $\aaa(f)$. Still his work also highlights how far the
current theoretical tools, for example \cite{Mi}, are from getting
close to the numerical conjectures. The same difficulty might also be
at work here.
\end{remark}

\begin{figure}[ht]
      \psfrag{p1}[][]{\begin{picture}(0,0)
            \put(-5, 0){\makebox(0,0)[l]{$z_1$}}
                      \end{picture}}
 \psfrag{p2}[][]{\begin{picture}(0,0)
            \put(-5, 0){\makebox(0,0)[l]{$z_2$}}
                      \end{picture}}
 \psfrag{p3}[][]{\begin{picture}(0,0)
            \put(-5, 0){\makebox(0,0)[l]{$z_3$}}
                      \end{picture}}
 \psfrag{p4}[][]{\begin{picture}(0,0)
            \put(-5, 0){\makebox(0,0)[l]{$z_4$}}
                      \end{picture}}
 \psfrag{p5}[][]{\begin{picture}(0,0)
            \put(-5,0){\makebox(0,0)[l]{$z_5$}}
                      \end{picture}}
 \psfrag{p6}[][]{\begin{picture}(0,0)
            \put(-5,0){\makebox(0,0)[l]{$z_6$}}
                      \end{picture}}
 \psfrag{p7}[][]{\begin{picture}(0,0)
            \put(-7, 0){\makebox(0,0)[l]{$z_7$}}
                      \end{picture}}
  \resizebox{0.8\linewidth}{!}{\includegraphics{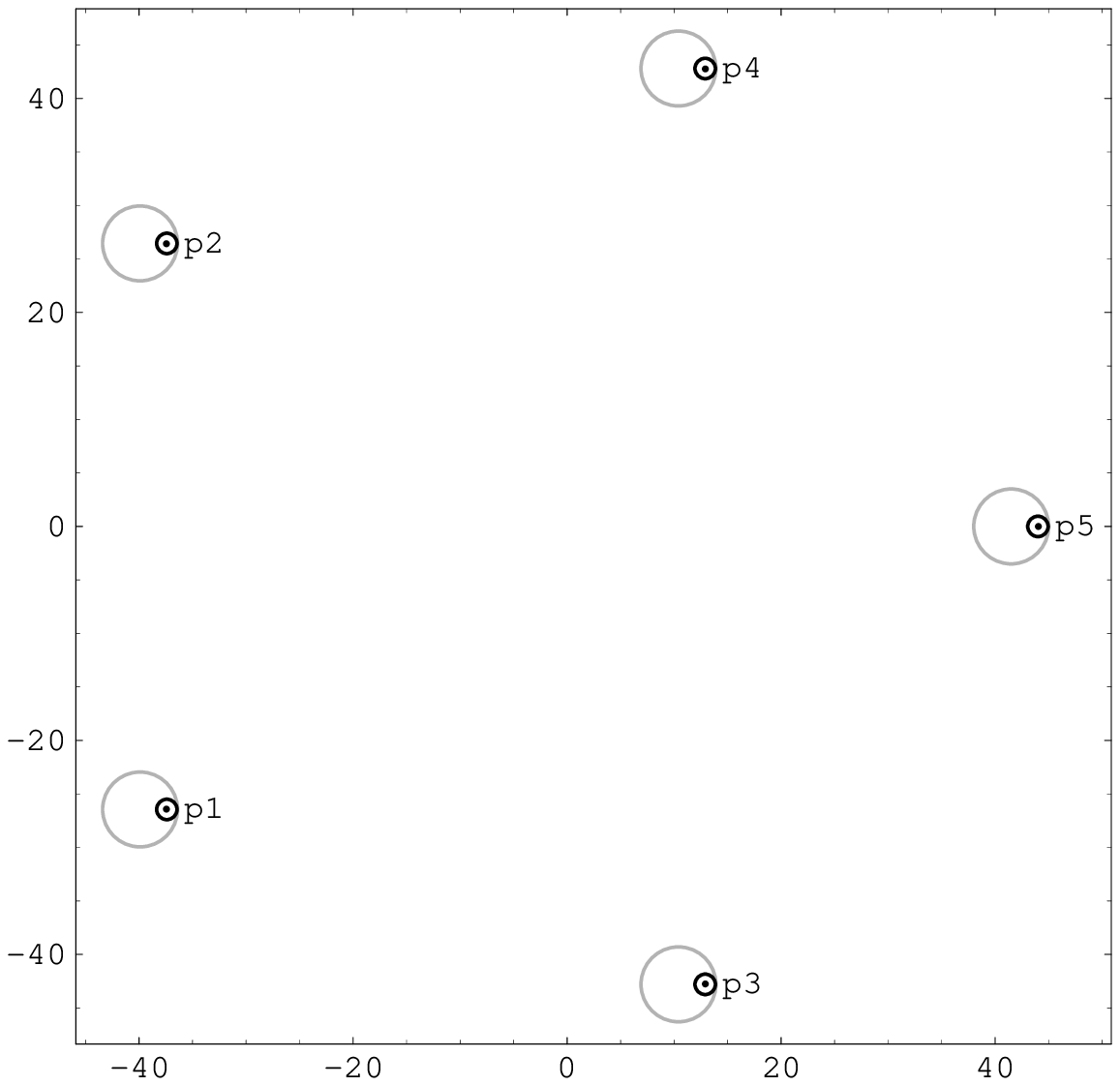}}
  \caption{A comparison of estimates for the roots of $z^5+5z^4-45^5$}

\medskip

      \psfrag{p1}[][]{\begin{picture}(0,0)
            \put(-97, 16){\makebox(0,0)[l]{$z_1$}}
                      \end{picture}}
  \resizebox{0.6\linewidth}{!}{\includegraphics{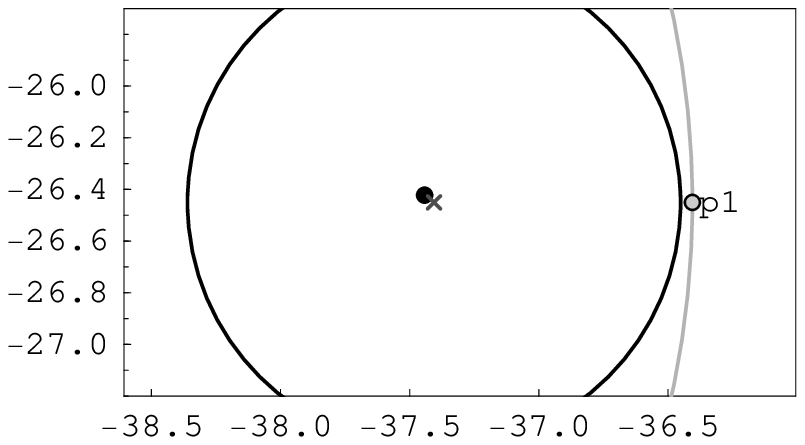}}
\caption{A zoom-in on $z_1$ with a root of $z^5-45^5$ marked by a
gray disk}
\end{figure}

\begin{example} \label{ex3}
Here we apply Corollary~\ref{cRS} to the polynomials  $\psi_{a,5}(z)
= z^5-a^5,\, a > 0$, and the operator $T=I + D$ to provide a precise  explanation for  the behavior of the roots exhibited in Figures~1 and~2. In this case we have
$\tau(\psi_{a,5}) = a$ and $(T\psi_{a,5})(z)=z^5+ 5 z^{n-1}-a^5$.
By Corollary~\ref{cRS}, if $a > 9$ we have
\begin{equation*} 
Z(T\psi_{a,5})
     \subset \{-1\} + Z(\psi_{a,5}) + \nD\bigl(\gamma_1 /a
\bigr),
\end{equation*}
where $\gamma_{1} = 42.944$. In Figures~1 and~2 we used $a=45$. Thus
\[
R_T(\psi_{45,5}) \leq \gamma_1/a = 0.95431\dot{1}.
\]
In terms of the geometric objects in Figures~1 and~2 the last
inequality predicts that the maximum distance between black dots and
crosses, paired in the natural way, is $< 0.954312$. In fact {\em
Mathematica} gives $0.046083$ as the maximum distance. Our estimate
therefore overestimates $0.046083$ by means of $0.954312$, which is
not impressive. But it can still be considered satisfactory because
the upper bound from Corollary~\ref{cRS}, or more generally from Theorem~\ref{lmt}, holds under the most general conditions.

Finally we notice that \cite[Corollary~5.4.1(iii)]{RS} in this
specific case would only give
\[
Z(T\psi_{45,5}) \subset \{-5/2\} + Z(\psi_{45,5}) + \nD\bigl( 5/2
\bigr),
\]
which is much less precise information. This is illustrated in
Figure~3, where the boundaries of the disks in the last inclusion are
gray, the smaller disks from our estimate are outlined in black, and
the roots of
\[
T\psi_{45,5} = z^5+5z^4-45^5
\]
are marked by black dots. In a magnified view in Figure~4 the gray
point marks $45 \exp(6 \pi i /5)$, while the cross is at $45 \exp(6
\pi i /5)-1$.
\end{example}

\medskip

\noindent{\bf Acknowledgment}. We are grateful to the referees
for many valuable suggestions that led to a much improved paper.

\end{document}